\documentclass[12pt,a4paper]{article}
\usepackage{amsmath, amssymb, theorem, latexsym}
\usepackage{graphicx}

\newtheorem{theorem}{Theorem}[section]
\newtheorem{lemma}[theorem]{Lemma}

\newtheorem{assumption}[theorem]{Assumption}

\newtheorem{corollary}[theorem]{Corollary}
\newtheorem{example}[theorem]{Example}
\newtheorem{remark}[theorem]{Remark}
\newtheorem{remarks}[theorem]{Remarks}
\newtheorem{property}[theorem]{Property}

\allowdisplaybreaks

\def \t{{\bf t}}
\def \v{{\bf v}}

\newcommand {\pa}{\partial}

\numberwithin{equation}{section}

\title{Subelliptic estimates for some systems of complex vector fields
: quasihomogeneous case}
\author{
M. Derridj\\
5 rue de la Juvini\`ere\\
78 350 Les loges en Josas, France.\\
and \\
B. Helffer \\
\upshape
Laboratoire de Math\'ematiques, Univ Paris-Sud and CNRS,\\
F 91 405 Orsay Cedex, France.}

\begin{document}
\maketitle
\newpage
\begin{abstract}
For about twenty five years it was a kind of folk theorem that complex
 vector-fields defined on $\Omega\times \mathbb R_t$ (with $\Omega$
 open set in $\mathbb R^n$) by
$$ L_j = \frac{\pa}{\pa t_j} + i \frac {\pa \varphi}{\pa t_j}(\t)\, \frac{\pa}{\pa
  x}\;,\;
 j=1,\dots, n\;,\; \t\in \Omega, x\in \mathbb R
\;,$$
with $\varphi$ analytic, 
were subelliptic as soon as  they were hypoelliptic.
This was the case when $n=1$ but in the case $n>1$,  an inaccurate
 reading of the proof given  by Maire (see also Tr\`eves) of the 
 hypoellipticity of such systems, under the condition that $\varphi$ does not admit
  any local maximum or minimum
 (through a non standard subelliptic estimate),  was supporting the
 belief for this folk theorem. Quite recently, J.L.~Journ\'e and
  J.M.Tr\'epreau
 show by examples that there are very simple systems (with polynomial
  $\varphi$'s) 
  which were hypoelliptic but not subelliptic in the standard $L^2$-sense. So it is natural
 to analyze this problem of subellipticity 
 which is in some sense intermediate (at least when $\varphi$ is $C^\infty$)
 between the maximal hypoellipticity
 (which was analyzed by Helffer-Nourrigat and Nourrigat)  and the
 simple local hypoellipticity (or local microhypoellipticity)  and to start first with the easiest non trivial  examples. 
The analysis presented here is a continuation of a previous work by
  the first author and is devoted to the case of quasihomogeneous functions.

\end{abstract}

 \newpage 
\section{Introduction and Main result}
\subsection{Preliminaries}
Let $\Omega$ an open set in $\mathbb R^n$ with $0\in \Omega$. 
We consider the regularity properties of the following family of complex vector fields
 on $\Omega \times \mathbb R$
\begin{equation}\label{sys}
L_j = \frac{\pa}{\pa t_j} + i \frac {\pa \varphi}{\pa t_j}(\t)\, \frac{\pa}{\pa
  x}\;,\;
 j=1,\dots, n\;,\; \t\in \Omega, x\in \mathbb R
\;,
\end{equation}
where $\varphi \in C^1(\Omega, \mathbb R)$, with $\varphi(0)=0$. We
 will concentrate our analysis
 near a point $(0,0)$ (but note that the operator is invariant by translation in the $x$ variable).\\

Many authors have considered this type of systems. For a given
$\Omega$, they were in
particular interested
 in the existence, for some pair  $(s,N)$
 such that $s+N >0$, of the following family of inequalities.\\ For any pair of bounded
 open sets $\omega$, $I$ such that $ \overline{ \omega} \subset
 \Omega$
 and  $I\subset \mathbb  R $, there exists a constant $ C_{s,N}(\omega, I)$
 such that
\begin{equation}\label{subells}
|| u||^2_s \leq C_N(\omega, I) \left( \sum_{j=1}^n ||L_j u||^2_0 + ||u||^2_{-N}\right)
\;,\; \forall u \in C_0^\infty (\omega\times I)\;,
\end{equation}
where $||\;\cdot\;||_r$ denotes the Sobolev norm in $H^r(\Omega \times
\mathbb R)$.\\
If $s>0$, we say that we have a subelliptic estimate. In \cite{Tr},
there are also results where $s$ can be arbitrarily negative. We will
speak in the case when $s\leq
0$ of  weak subellipticity. Note that in this case  the existence of this inequality
 is not sufficient for proving hypoellipticity.

The system \eqref{sys} being elliptic in the $\t$ variable, it is enough
 to analyze the subellipticity microlocally near $\tau =0$,
 i.e. near $(0, (0,\xi))$ in $ (\omega \times I)\times (\mathbb R^{n+1}\setminus \{0\})$
 with $\{\xi >0\}$  or $\{\xi <0\}$.\\
This leads to the analysis of the existence of two constants
 $C_s^+$ and $C_s^-$ such that the  two following inequalities
 hold~:
\begin{equation} \label{I+}
\int_{\omega \times \mathbb R^+} \xi^{2s} |\widehat u (\t,\xi)|^2
\, dt d\xi \leq C_{s}^+ \int_{\omega \times \mathbb R^+} | \widehat
   {Lu} (\t,\xi)|^2\; d\t d\xi\;,\; \forall u \in C_0^\infty
   (\omega\times \mathbb R)\;,
\end{equation}
where $\widehat u (\t,\xi)$ is the partial Fourier transform of $u$
 with respect to the $x$ variable, 
 and 
\begin{equation} \label{I-}
\int_{\omega \times \mathbb R^-} |\xi|^{2s} |\widehat u (\t,\xi)|^2
\, d\t d\xi \leq C_{s}^- \int_{\omega \times \mathbb R^-} | \widehat
   {Lu} (\t,\xi)|^2\; d\t d\xi\;,\; \forall u \in C_0^\infty
   (\omega\times \mathbb R)\;.
\end{equation}

When \eqref{I+} is satisfied, we will speak of microlocal
subellipticity
 in $\{\xi >0\}$ and similarly when \eqref{I-} is satisfied, we will
 speak of microlocal subellipticity in $\{\xi <0\}$. Of course, when $s
 >0$, it is standard that these two inequalities imply
 \eqref{subells}.\\
We now  observe that \eqref{I+} for $\varphi$
 is equivalent to \eqref{I-} for $-\varphi$, so it is enough to
 concentrate our analysis on the first case.

\subsection{The main result}
In \cite{De}, the first author gave a sufficient condition on
$\varphi$
 for getting \eqref{subells} with $s>0$. In this article, we consider
 the case of quasihomogeneous functions $\varphi$ on $\mathbb R^2$
 (i.e. $n=2$) and we will give a simple condition of subellipticity 
 where $s$ will be related rather simply with the quasihomogeneity of
 $\varphi$.

These conditions will be expressed for $\varphi$ in $C^1$ but note that
they become more simple in the analytic case (see Section \ref{s6}).

More precisely, let  $m$ and $\ell$ in
$\mathbb R^+$ such that
\begin{equation}\label{condl}
\ell \geq 1 \;, 
\end{equation}
and
\begin{equation}\label{condlm}
m\geq 2 \ell \;.
\end{equation}
{\bf We make these two assumptions\footnote{
with in addition the assumption that $\ell$ is rational in the
analytic case,} in the whole paper.}
 
We consider in $\mathbb R^2$ $(t,s)$ as the variables (instead of
$\t$) and  the functions $\varphi\in
C^1(\mathbb R^2)$ will be quasihomogeneous in the following sense
\begin{equation}\label{quash}
\varphi (\lambda t,\lambda ^\ell s) = \lambda^m \varphi (t,s)\;,\;
 \forall (t,s,\lambda)\in \mathbb R^2\times \mathbb R^+\;.
\end{equation}

According to \eqref{quash}, the real function $\varphi$ is determined
by its restriction $\widetilde \varphi$ to the distorted circle
$\mathcal S$
\begin{equation}\label{defphitilde}
\widetilde \varphi := \varphi_{\vert \mathcal  S}\;.
\end{equation}
where $\mathcal S$ is defined by
\begin{equation}\label{defS}
\mathcal  S = \{(t,s)\,;\, t^{2\ell}+ s^2 =1\}\;,
\end{equation}

Our main result is stated under the following assumption on
$\widetilde \varphi$.\\

\begin{assumption}{ (H2)} \label{H2}
\begin{enumerate}
\item $\widetilde \varphi$ is not strictly negative.
\item $\widetilde \varphi$ cannot have, at any of its zeroes,  a local maximum.
\item If $\mathcal S_j^+$ is a  connected component of $\widetilde
  \varphi^{(-1)} (]0,+\infty[)$, then one can write $\mathcal S_j^+$
as a finite union of arcs satisfying  Property \ref{Prop1} below.
\item If $\mathcal S_j^-$ is one connected component of $\widetilde
  \varphi^{(-1)} (]-\infty,0[)$, then $\tilde \varphi$ has a unique
  minimum in $\mathcal S_j^-$.
\item There exists $p\geq 1$, such that, if $\theta_0$ is a zero of
  $\widetilde \varphi$, then there exists an open arc  $\mathcal
  V_{\theta_0}$
 containing  $\theta_0$ in $\mathcal S$
  and $C_0 >0$, such that
\begin{equation}\label{4.1}
| \widetilde \varphi (\theta) - \widetilde \varphi (\theta')|\geq \frac{1}{C_0}
 |\theta - \theta'|^p\;,\; \forall \theta, \theta'\in \mathcal
 V_{\theta_0}\;,
\end{equation}
with $\theta$ and $\theta'$ in the same connected component in $\mathcal
V_{\theta_0}\setminus \{\theta_0\}$.
\end{enumerate}
\end{assumption}

Here in the third item of $(H2)$, we mean by saying that a closed arc
$\left[\theta,\theta'\right]$ has Property \ref{Prop1}  the following~:
\begin{property}\label{Prop1}~\\
There exists on this arc  a point $\widehat \theta$ such that~:
\begin{enumerate}
\item[(a)]
$\widetilde \varphi$ is non decreasing on the arc  $\left[\theta, \widehat \theta\right]$
 and  non increasing on the arc $\left[\widehat \theta
   ,\theta'\right]$. (So the restriction of 
 $\widetilde{ \varphi}$ to  $\left[\theta,\theta'\right]$  has a maximum at $\widehat \theta$).
\item[(b)] \begin{equation}\label{ineglarge}
 \langle \widehat \theta\,|\,
 \theta  \rangle_\ell  \geq 0\; \mbox{and}\;   \langle \widehat \theta\,|\,
 \theta'  \rangle_\ell  \geq 0\;,        \end{equation}
where for $\theta =(\alpha,\beta)$ and  $\widehat \theta =(\widehat
 \alpha,\widehat\beta)$ in $\mathcal S\subset \mathbb R^2$, 
\begin{equation}
\langle \widehat \theta\,|\,
 \theta  \rangle_\ell :=\widehat \alpha \alpha  |\widehat \alpha|^{\ell
 -1} |\alpha|^{\ell-1}
 + \widehat \beta \beta\;.
\end{equation} 
\end{enumerate}
\end{property}
Note here that we could have $\widetilde \varphi$ constant on
$\mathcal S_{\theta,\theta'}$ and $\widehat \theta = \theta$ or
$\theta'$. Moreover item (b)  says that
the length of the two arcs is sufficiently small, more precisely that
  the distorted
``angles''
 (see Section \ref{s2}) associated to $[\theta,\theta']$
 are acute.\\
We can now state our main theorem~:
\begin{theorem}\label{theogen}~\\
Let $\varphi \in C^1(\mathbb R^2,\mathbb R)$ satisfying
\eqref{quash}, with $\ell$ and $m$ satisfying \eqref{condl} and
\eqref{condlm}. Then $L$ is microlocally subelliptic 
in $\{\xi >0\}$.
\end{theorem}

\begin{remarks}~
\begin{enumerate}
\item The proof of Theorem \ref{theogen} consists
 in showing that  Assumption $(H2)$  implies the  assumption $(H_+(\alpha))$
 with $\alpha = \frac {1}{\max(m,p)}$ introduced in \cite{De}
 and which will be recalled in Section \ref{DeCr}.
\item 
 \cite{De} was considering  the homogeneous case $\ell
 =1$
 and $m\geq 2$. Here we generalize by considering the quasihomogeneous
  case but the sufficient condition given here for getting 
 Assumption \ref{Ass1} satisfied will be already an improvement in the
 homogeneous case.
\item The conditions on $\widetilde \varphi$ are more restrictive on
 the
 arcs $S_k^-$.
\item If $\varphi$ is analytic and $\ell$ is rational, 
 the statement of the main theorem becomes simpler. (iii) and (v)
 are indeed automatically satisfied as soon that $\widetilde \varphi$
 is not identically $0$. Moreover, if we write 
$\ell= \frac{\ell_2}{\ell_1}$ (with $\ell_1$ and $\ell_2$ mutually
 prime integers),
 all the criteria involving $\widetilde \varphi$ can be reinterpreted
 as criteria involving  the restriction $\widehat \varphi$
 of $\varphi$
 on $$
\mathcal S_{\ell_1,\ell_2} =\{ (t,s)\;;\; t^{2\ell_2} + s^{2\ell_1}=1\}\;.
$$
\item
The condition $(H2), ii)$ 
 is natural and cannot  be relaxed
 according to the necessary conditions obtained by  J.L. Journ\'e and
 J.M.~Tr\'epreau \cite{Tr}   for the subellipticity of these systems.

\end{enumerate}
\end{remarks}

\paragraph{Organization of the paper}~\\
The proof of our main theorem  will be based on a ``abstract'' criterion established in
\cite{De},
 which will be recalled in Section \ref{DeCr}. After the introduction
 of a terminology adapted to the quasihomogeneity of the problem
 in Section \ref{s2}, we continue with the proof of the general
 main theorem in three steps starting from the analysis of  the quasielliptic case in
 Section \ref{s3},  showing then how
 one can localize the proof in
 suitable quasihomogeneous sectors in Section \ref{s5a} and finishing by proving the
 general case in Section \ref{s5}.
Section \ref{s6} is devoted to the particular case of an analytic
function $\varphi$. We give in an appendix
 the computation of a basic Jacobian, whose control is important
 in the verification of the abstract criterion.

{\bf Acknowledgements}\\
At various stages of this work, the authors have fruitful discussions
with
 F.~Nier, 
 J.M.~Tr\'epreau and H.~Maire.

\section{Derridj's subellipticity criterion.}\label{DeCr}
\subsection{The statement}
We now recall the criterion established in \cite{De}. This involves,
for a given $\alpha >0$,  the following geometric escape condition on
$\varphi$. We do not have in this section the restriction $n=2$

\begin{assumption}{$(H_+(\alpha))$}\label{Ass1}~\\
There exist an open set $\omega\subset \Omega$ and  $\widetilde \omega \subset \omega$, with full Lebesgue
measure in $\omega$, and a map
 $\gamma$~: 
$$
\widetilde \omega\times [0,1]\ni (\t,\tau)\mapsto
\gamma(\t,\tau) \in \Omega\;,
$$
 such that
\begin{enumerate}
\item $\gamma (\t,0)=\t\,;\, \gamma (\t,1) \not\in \omega\;,\, \forall
\t\in \widetilde \omega\;$.
\item $\gamma$ is of class $C^1$ outside a negligeable set $E$ and
  there exist  $C_1>0$, $C_2>0$  and  $C_3>0$ such that
\begin{enumerate}
\item $$
|\pa_\tau \gamma (\t,\tau)|\leq C_2\;,\; \forall (\t,\tau) \in \widetilde
  \omega \times [0,1]\setminus E\;.
$$
\item
\begin{equation}\label{minorbas}
|\det (D_\t \gamma)(\t,\tau)|\geq \frac {1}{ C_1}\;,
\end{equation}
where $D_\t \gamma$ denotes the Jacobian matrix of $\gamma$
 considered as a map from $\widetilde \omega$ into $\mathbb R^2$.

\item
$$ \varphi(\gamma(\t,\tau)) -\varphi(\t) \geq \frac{1}{ C_3} \tau^\alpha\;,\;
  \forall (\t,\tau) \in \widetilde \omega \times [0,1]\;.
$$
\end{enumerate}
\end{enumerate}
\end{assumption}

Using this assumption, it is proved in  \cite{De} the following theorem.
\begin{theorem}\label{thde}~\\
If $\varphi$ satisfies $(H_+(\alpha))$, then the 
 associated system \eqref{sys}$_\varphi$ is microlocally \break
$\frac{1}{\alpha}$-subelliptic
 in $\{\xi >0\}$.\\
If $-\varphi$ satisfies $(H_+(\alpha))$, then the 
 associated system \eqref{sys}$_\varphi$ is microlocally \break 
$\frac{1}{\alpha}$-subelliptic
 in $\{\xi <0\}$.\\
If $\varphi$ and $-\varphi$ satisfy  $(H_+(\alpha))$, then the 
 associated system \eqref{sys}$_\varphi$ is \break 
$\frac{1}{\alpha}$-subelliptic.\end{theorem}

\subsection{The proof}
For the commodity of the reader, we reproduce the proof of \cite{De}
 in the case $\xi >0$. \\
If  $u\in C_0^\infty(\omega \times \mathbb R)$, one can,
 using the partial Fourier transform (with respect to
 $x$), recover $u$ from $f=L u$ by 
\begin{equation}\label{sol}
\widehat u (\t,\xi) = -\int_0^1 \exp \left[ \xi 
  \varphi(\gamma(\t, \tau) \right] \;
 \widehat f ( \gamma(\t, \tau),\xi)\cdot (\pa_\tau \gamma) (\t,\tau)\; d\tau\;.
\end{equation}
Taking Cauchy-Schwarz in \eqref{sol}, we obtain
\begin{equation}\label{casc}
\begin{array}{l}
|\widehat u (\t,\xi)|^2  \leq\\
\quad  \left(\int_0^1 \exp \left[ \xi \,
  \varphi(\gamma(\t, \tau))\right] d\tau \right) \; \times\\
\qquad\qquad \times 
\left(\int_0^1 \exp \left[ \xi \,
  \varphi(\gamma(\t, \tau))\right] \,
 |\widehat f ( \gamma(\t, \tau),\xi)\cdot (\pa_\tau \gamma) (\t,\tau)\;|^2
 d\tau\right)\;.
\end{array}
\end{equation}
By  items   $(ii)_{(b)}$ and  $(ii)_{(c)}$ in Assumption \ref{Ass1}, this implies
$$
|\widehat u (\t,\xi)|^2  \leq C_2^{\,2}
\quad  \left(\int_0^1 \exp -\frac{1}{C_3}  \tau^\alpha \xi \;d \tau \right) \; \times\\
\left(\int_0^1 \exp - \frac{1}{C_3}  \tau^\alpha \xi \;\,
 |\widehat f ( \gamma(\t, \tau),\xi)|^2\,
 d\tau\right)\;.
$$
So, integrating in $\t$ over   $\omega$, we get
\begin{equation}\label{23a}
\begin{array}{l}
\int_\omega |\widehat u (\t,\xi)|^2 \;d\t \\
\quad  \leq C_2^2
\quad  \left(\int_0^1 \exp -\frac{1}{C_3} 
 \tau^\alpha \xi \;d \tau \right) \; \times
\left(\int_0^1\int_\omega  \exp -\frac{1}{C_3}  \tau^\alpha \xi \;
 |\widehat f ( \gamma(\t, \tau),\xi)|^2
 d\t d\tau\right)\;.
\end{array}
\end{equation}
We now change of variables~: $\v=\gamma(\t,\tau)\,,\, \tau = \tau$.
The second term in the r.h.s. of \eqref{23a} can be estimated as follows.
$$
\begin{array}{l}
\int_0^1\int_\omega  (\exp -\frac{1}{C_3} \tau^\alpha \xi) \;
 |\widehat f ( \gamma(\t, \tau)),\xi)\;|^2
 d\t d\tau \\
\quad\quad  = \int_0^1\int_\omega  (\exp -\frac{1}{C_3} \tau^\alpha \xi) \;
 |\widehat f (\v,\xi)\;|^2 |D_\t\gamma|^{-1}
 d\tau d\v\\
\quad\quad \leq C_1  \left(\int_0^1  \exp -\frac{1}{C_3} \tau^\alpha
\xi\,d\tau \right)
 \; \left( \int_\Omega  |\widehat f (\v,\xi)|^2 d\v\right)\;,
\end{array}
$$
where we have used the lower bound  \eqref{minorbas} for the Jacobian 
$|D_\t\gamma|$\;.

So finally, one has 
$$
\int_\omega |\widehat u (\t,\xi)|^2 \;d\t  \leq C_1 C_2^2  
\; \left(\int_0^1 (\exp -\frac{1}{C_3} \tau^\alpha \xi) \,d \tau
\right)^2 \left( \int_\Omega  |\widehat f (\v,\xi)|^2 d\v\right)\;.
$$

We then obtain the existence of $C (C_1,C_2,C_3) >0$ such that, for all $\xi >0$,
$$
\int_\omega |\widehat u (\t,\xi)|^2 \;d\t  \leq C (C_1,C_2,C_3)  |\xi|^{-\frac 2
  \alpha}\;\left( \int_\Omega  |\widehat f (\v,\xi)|^2 d\v\right)\;.
$$

\section{Quasihomogeneous structure}\label{s2}
\subsection{Distorted geometry}
Condition (i) in Assumption  \ref{Ass1} expresses the property that
the curve is escaping from $\omega$. 
For  the description of escaping curves, it appears useful to 
extend  the usual terminology used in the Euclidean space $\mathbb
R^2$
 in a way which is adapted to the given quasihomogeneous structure.
 This is realized by introducing  the {\it dressing} map~:
\begin{equation}\label{dressing}
(t,s)\mapsto d_\ell (t,s)= \left( t\,|t|^{\ell -1},s\right)\;.
\end{equation}
which is at least of class $C^1$ as $\ell \geq 1$, and 
 whose main role is to transport the distorted geometry onto the
 Euclidean
 geometry.\\
The first example was the unit distorted circle (in short unit
 disto-circle
 or unit ``circle'') $\mathcal S$ introduced in \eqref{defS}
 whose image by $d_\ell$ becomes the standard unit circle in $\mathbb
 R^2$
 centered at $(0,0)$.\\
Similarly, we will speak of
disto-sectors, disto-arcs, disto-rays, disto-disks. In particular, for $(a,b)\in \mathcal S$, we
define
 the disto-ray $\mathcal R_{(a,b)}$ by
\begin{equation}
\mathcal R_{(a,b)} :=\{ (\lambda a,\lambda^\ell b)\,;\, 0 \leq \lambda \leq
1\}\;.
\end{equation}

The disto-scalar product of two vectors in $\mathbb R^2$
 $(t,s)$ et $(t',s')$ is then given by
\begin{equation} \label{defdist}
\langle (t,s)\;|\; (t',s')\rangle_{\ell} = tt' |tt'|^{\ell -1} +
ss'\;.
\end{equation}
(for $\ell=1$, we recover the standard scalar product).\\
For $(t,s)\in \mathbb R^2$, we introduce also the quasihomogeneous
positive
 function $\varrho$ defined on $\mathbb R^2$ by~:
\begin{equation}\label{defrho}
\varrho (t,s)^{2\ell} = t^{2\ell}  + s^2\;.
\end{equation}
With these notations, we observe that, if $(t,s)\in \mathbb
R^2\setminus \{(0,0)\}$, then
\begin{equation}
(\widetilde t\,,\, \widetilde s):= \left( \frac{t}{\varrho(t,s)}\,,\, \frac{s}{\varrho(t,s)}\right)\in
\mathcal S\;,
\end{equation}
and
$$
(t,s) \in \mathcal R_{(\widetilde t, \widetilde s)}\;.
$$
The open disto-disk $D(R)$ is then defined by
\begin{equation}\label{defdisto}
D(R)=\{(x,y)\;\vert\; \varrho(x,y) < R\}\;.
\end{equation}
We can also consider
a parametrization of 
  the
disto-circle by a parameter on the corresponding 
circle $\vartheta \in \mathbb R
/2\pi \mathbb Z$ (through the dressing map). We note
 that we have a 
 natural (anticlockwise) orientation of the disto-circle. In other cases it will be better to parametrize by $s$ (if $t\neq 0$) or by $t$ (if $s\neq 0$). So a point in $\mathcal S$ will be defined either by $\theta$
 or by $(a,b)\in \mathbb R^2$ or by $\vartheta$.\\

Once an orientation is defined on $\mathcal S$, two points $\theta_1$ and $\theta_2$ (or $(a_1,b_1)$ and
 $(a_2,b_2)$) on $\mathcal S$ will determine a unique unit ``sector'' $V\subset
 D(1)$. \\

\subsection{Distorted dynamics}
The parametrized curves $\gamma$ permitting us to satisfy Assumption 
\ref{Ass1}
 will actually be ``lines'' (possibly broken) which will finally
 escape
 from a neighborhood of the origin. Our aim in this subsection is to
 define  these ``lines'' (actually
 distorted parametrized lines).\\
In parametric coordinates, with  
\begin{equation}\label{deft} t(\tau)=t+\varrho \tau\,,
\end{equation} the curve $\gamma$
 starting from $(t,s)$ and disto-parallel to $(c,d)$ is defined by
 writing that the vector 
$(t(\tau)|t(\tau)|^{\ell -1}- t|t|^{\ell -1} , s(\tau)-s)$ is parallel
 to  $(c |c|^{\ell -1}, d)$.\\
{\bf In the applications, we will consider only consider $\varrho =\pm c$.}\\
So
$$
\left(t(\tau)|t(\tau)|^{\ell -1}- t|t| ^{\ell -1}\right) d = c
|c|^{\ell -1} (s(\tau)-s)\;,
$$
and we find 
\begin{equation}\label{defs}
s(\tau) = s + \frac {d}{c|c|^{\ell -1}} \left( t(\tau)|t(\tau)|^{\ell
    -1}- t| t| ^{\ell -1}\right)\;,
\end{equation}
We consider the map $\sigma \mapsto f_\ell (\sigma)$ which is  defined by
\begin{equation}
f_\ell(\sigma) = \sigma |\sigma|^{\ell -1}\;.
\end{equation}
Note that
\begin{equation}
f_\ell'(\sigma) = \ell |\sigma|^{\ell -1}\geq 0\;.
\end{equation}
With this new function, \eqref{defs} can be written as
\begin{equation}\label{inva}
d f_\ell(t(\tau)) - s(\tau) f_\ell (c) = d f_\ell(t) - s f_\ell (c)\;
.\end{equation}

This leads us to use the notion of distorted determinant
 of two vectors in $\mathbb R^2$. For two vectors $v:=(v_1,v_2)$
 and $w:=(w_1,w_2)$, we introduce~:
\begin{equation}\label{defdet}
\Delta_\ell (v;w) = f_\ell(v_1) w_2 - v_2 f_\ell(w_1)\;.
\end{equation}
We will also write~:
\begin{equation}\label{defdetbis}
\Delta_\ell (v;w) =\Delta_\ell(v_1,v_2,w_1,w_2)\;.
\end{equation}
With these notations, \eqref{inva} can be written
\begin{equation}
\Delta_\ell((c,d);(t(\tau),s(\tau))) = \Delta_\ell ((c,d);(t,s))\;,
\end{equation}
or
\begin{equation}
\Delta_\ell(c,d,t(\tau),s(\tau)) = \Delta_\ell (c,d,t,s)\;,
\end{equation}

When $\ell=1$, we recover the usual determinant of two vectors
 in $\mathbb R^2$. For $\ell \geq 1$, we have simply the relation~:
\begin{equation}
\Delta_\ell (v;w) = \Delta_1( d_\ell(v);d_\ell(w))\;.
\end{equation}
Note that $\Delta_\ell (v;w)$ vanishes
 when $d_\ell(v)$ and $d_\ell(w)$ are collinear.\\

We now look at the variation of $\psi$ which 
 is defined (for a given initial point $(t,s)$) by 
\begin{equation}
\tau \mapsto \psi(\tau) = \rho(\tau)^{2\ell} = t(\tau)^{2\ell} + s(\tau)^2\;.
\end{equation} 
Easy computations give also~:
\begin{equation}
\psi'(\tau) =\frac{2\varrho}{f_\ell(c)} f'_\ell (t+ \varrho  \tau) \langle (c,d) \;|\;
(t(\tau), s(\tau))\rangle_{\ell}\;,
\end{equation}
whose sign is the product of the sign of the (disto)-scalar product
 of $(c,d)$ and $(t(\tau),s(\tau))$ and  the sign of
 $c\varrho$.\\

We now analyze the variation of the (disto)-scalar product
$\langle (c,d) \;|\;
(t(\tau), s(\tau))\rangle_{\ell}$ as a function of $\tau$.
We have the formula
\begin{equation}\label{varprosca}
\langle (c,d) \;|\;
(t(\tau), s(\tau))\rangle_{\ell} = \langle (c,d) \;|\;
(t, s)\rangle_{\ell} + \frac{1}{f_\ell (c)}(f_\ell (t(\tau) )- f_\ell (t))\;.
\end{equation}
If we now assume that 
\begin{equation}\label{hip1}
c\varrho >0\;,\;  \langle (c,d) \;|\;
(a, b)\rangle_{\ell} \geq 0 \;,
\end{equation}
Then for $(s,t)$ in the unit sector $\mathcal V_{abcd}$ associated to the arc
$\left((a,b)\,,\, (c,d)\right)$, we obtain~:
$$
\psi'(\tau) \geq \frac{1} {f_\ell(c)^2} \,\times
\, (2\varrho f'_\ell(t+\varrho \tau))\,\times\, (f_\ell (t(\tau) )- f_\ell
 (t))\;.
$$
We rewrite this inequality in the form
$$
\psi'(\sigma) \geq \frac{1} {f_\ell(c)^2} \,\times
 \left( (f_\ell (t(\sigma) )- f_\ell
 (t))^2\right)'\;,\;\forall \sigma \geq 0.
$$
Integrating over $[0,\tau]$, we get for $\tau \geq 0$~:
\begin{equation}
\psi(\tau) \geq \frac{1} {f_\ell(c)^2} \,\times
  (f_\ell (t(\tau) )- f_\ell
 (t))^2\;.
\end{equation}
We now need the following
\begin{lemma}\label{LF1}~\\
For any  $\ell \geq 1$,  $ \tau \geq 0$, and $
\gamma \in \mathbb R$, we have
\begin{equation}\label{inequt3}
f_\ell(\tau + \gamma) - f_\ell (\gamma) \geq  f_\ell (\frac \tau 2)\;.
\end{equation}
\end{lemma}
\begin{remark}~\\
This lemma can be improved  when $\gamma \geq 0$; we can then show
\begin{equation}\label{inequt4}
f_\ell(\tau + \gamma) - f_\ell (\gamma) \geq  f_\ell (\tau)\;.
\end{equation}
\end{remark}
{\bf Proof}~\\
By the previous remark, the proof is clear
 when $\gamma \geq 0$ or  when $\gamma +\tau \leq 0$.
So it remains to analyze the case when
$  - \tau <\gamma <0$. But the two terms on the left hand side are now 
 positive. So we immediately obtain \eqref{inequt3} 
by observing that $\max ( \tau + \gamma, - \gamma) \geq \frac \tau 2$.

\begin{remark}~\\
If we take the square, we obtain (and this time for any $\tau\in
\mathbb R$) the inequality
\begin{equation}
\left(f_\ell(\tau + \gamma) - f_\ell (\gamma)\right)^2
\geq  (\frac \tau 2)^{2\ell}\;.
\end{equation}
\end{remark}
Now using Lemma \ref{LF1}, this leads to
\begin{lemma}~\\
Under Condition \eqref{hip1}, we have, for any $\tau \geq 0$, for any
$(t,s)\in \mathcal V_{abcd}$, 
\begin{equation}\label{inequt1}
\rho(\tau)^{2 \ell} - \rho(0)^{2\ell}  \geq 
(\frac{\varrho \tau} {2c})^{2 \ell}\;.
\end{equation}
If instead $\varrho c <0$, we obtain~:
\begin{equation}\label{inequt2}
\rho(\tau)^{2 \ell} - \rho(0)^{2\ell} \leq - 
(\frac{\varrho \tau} {2c})^{2 \ell}\;.
\end{equation}
\end{lemma}

We continue by analyzing the variation of $s(\tau)$ and $t(\tau)$
 and more precisely the variation on the disto-circle of~:
\begin{equation}
\widetilde t (\tau) = \frac{t(\tau)}{\rho(\tau)}\;,\;
 \widetilde s(\tau) = \frac{s(\tau)} {\rho(\tau)^\ell}\;.
\end{equation}
After some computations, we get, with
 $$ \varrho = \pm c\;,$$
\begin{equation} \label{ttildeprimea}
\widetilde t' (\tau) = \pm |c|^{1-\ell} \frac{s(\tau)}{\rho(\tau)^{2\ell
    +1}} \Delta_\ell(c,d,t,s)\;,
\end{equation}
which can also be written in the form
\begin{equation} \label{ttildeprimeb}
\widetilde t' (\tau) = \pm |c|^{1-\ell} \frac{\widetilde s(\tau)}{\rho(\tau)}
\Delta_\ell (c,d,\widetilde t(\tau), \widetilde s(\tau))\;.
\end{equation}
Similarly, we get for $\widetilde s'$, 
\begin{equation} \label{stildeprimea}
\widetilde s' (\tau) = \mp \ell  |c|^{1-\ell} \frac{t(\tau)^{2\ell
    -1}}{\rho(\tau)^{3\ell}}\Delta_\ell (c,d,t,s)\;,
\end{equation}
and
\begin{equation} \label{stildeprimeb}
\widetilde s' (\tau) = \mp  \ell |c|^{1-\ell} \frac{\widetilde
  t^{2\ell -1} (\tau)}{\rho(\tau)}\Delta_\ell (c,d, \widetilde t (\tau),\widetilde s (\tau))\;.
\end{equation}
\section{Analysis of the quasielliptic 
 case ($\widetilde \varphi >0$)} \label{s3}
We first start the proof of the main theorem with the particular case when
\begin{equation}\label{strictpos}
\widetilde \varphi \geq \mu >0\;.
\end{equation}
This case is already interesting for presenting the main ingredients
of the general proof. We can remark  indeed that what we are doing below
 in $\mathcal S$ can be done later  in a specific (disto)-arc
 of $\mathcal S$.

On $\mathcal S$, taking a regular parametrization of $\mathcal S$
denoted
 by $\theta$, we consider the connected components of $\{ \widetilde
\varphi' \geq 0\}$ or of $\{\widetilde
\varphi' \leq 0\}$.

Our assumption takes in this case the following form.
\begin{assumption}{$(H_1)$}~\\
Either the cardinal of the connected components in $\mathcal S$
 of $\{\widetilde
\varphi' \geq 0\}$ or the cardinal of the  connected components of $\{\widetilde
\varphi' \leq 0\}$ is finite.
\end{assumption}
\begin{remark}\label{reman}~\\
This assumption is automatically satisfied if $\widetilde \varphi$ is
analytic.
\end{remark}

\begin{theorem}\label{thqh}~\\
Let $\varphi$ in $C^1(\mathbb R^2;\mathbb R)$ satisfying \eqref{quash} and \eqref{strictpos}.
Then $(H_1)$ implies $(H_+(\alpha))$, with $\alpha=\frac 1m$.
\end{theorem}
\begin{corollary}\label{th3.1}~\\
If $(H_1)$ is satisfied, then the system \eqref{sys}
 is $\frac 1m$-subelliptic in $\{\xi >0\}$\,.
\end{corollary}
~\\
{\bf Proof of Theorem \ref{thqh}}

\paragraph{Step 1 : Construction of the covering}~\\
Under Assumption $(H1)$ (say for definiteness under the first alternative),
 and considering the connected components introduced there,  we first start by constructing  a  finite
 covering
 of $\mathcal S$ by a family of disjoints open arcs $\mathcal S_j$
 such that
\begin{equation}\label{cov}~\\
\mathcal S = \cup_j \overline{\mathcal S_j}\;,
\end{equation}
in the following way. 
If $\theta_j$ (or $(a_j,b_j)$)  denotes the sequence of the left end points
 of the components of $\{ \widetilde
\varphi' \geq 0\}$ and by  $\widehat \theta_j$ (or $(c_j,d_j)$) the
 sequence of the right  end points, we define  $\mathcal S_j$ 
as the arc 
 $
\mathcal S_j = (\theta_j,\theta_{j+1})$.  We observe  that~:\\
{\it $\widetilde \varphi$ is non decreasing on $(\theta_j, \widehat \theta_j)$
 and strictly  decreasing on $(\widehat \theta _j,\theta_{j+1})$.}\\
We  now associate to the (disto)-arc $\mathcal S_j$ its (disto)-unit
sector $\mathcal V_j$.
For technical reasons, we will add a finite number  of points 
 such that we finally a (possibly new) family
 of open arcs $\mathcal S_j= (\theta_j,\theta_{j+1})$ such that each
 arc satisfies Property~\ref{Prop1}.

\paragraph{Step 2 : Construction of $\gamma$}~\\
We construct $\gamma$ independently in each  sector $\mathcal V_j$. More precisely, this will
 be a map from 
 $\left[(\mathcal V_j \cup \mathcal R_{\theta_j})\setminus \{(0,0)\}\right]\times
 [0,1]$ into $\Omega$ (and actually in the infinite sector associated to $\mathcal S_j$). From now on in this paragraph,  we fixed $j$ (and take it equal to  $1$). 
So $\mathcal S_1$ denotes the set of points $(t,s)\in D(1)$ such that
  $(\frac{ t} {\rho(t,s)}, \frac{ s}  {\rho(t,s)^\ell})\in
 (\theta_1,\theta_2)$.

We now define $\gamma$ (see \eqref{deft}-\eqref{defs}, with
$\varrho=c_1$) by
\begin{enumerate}
\item If $c_1 \neq 0$,
\begin{equation}\label{304a} \gamma (t,s,\tau):= (t(\tau),s(\tau)):= 
(t +  c_1 \tau, s + \frac{d_1}{f_\ell(c_1)} (f_\ell(t(\tau)) - f_\ell(t))\;.
\end{equation}
\item If $c_1 = 0$
\begin{equation} \label{304b}
 \gamma (t,s,\tau):= (t(\tau),s(\tau)):=(t, s + d_1 \tau)\;.
\end{equation}
\end{enumerate}
\begin{remarks}~\\
\begin{enumerate}
\item
Note that for any $(t_0,s_0,\tau)$ the
Jacobian
 of the map $(t,s)\mapsto \gamma(t,s,\tau)$ at $(t_0,s_0)$ is $1$.
\item
  Actually, one can avoid the second case by replacing $c_1=0$ by
 an arbitrarily close $\widetilde c_1$ whose sign will depend
 on the considered ``sector''.
\end{enumerate}
\end{remarks}
~From now on, we write for simplification 
$c=c_1, d= d_1$.\\

Let us look at the most generic case\footnote{We will  add in
  footnotes some  indications for the other case.} when $c\neq 0$.
 In order to show $(H_+(\alpha))$, the only non trivial property
is to show property (ii) (b)  in Assumption \ref{Ass1}.
\\
\begin{lemma}\label{lemmastays}~\\
If $(t,s)$ belongs to the ``subsector''  associated with $[\theta_1,\widehat
\theta_1]$
(resp. to the ``subsector''  associated with $[\widehat \theta_1,
\theta_2]$), then  the whole curve $\gamma (t,s,\tau)$ stays in the same
(infinite) ``subsector''
 for $\tau \geq 0$.
\end{lemma}
{\bf Proof}~~\\
The lemma is geometrically evident using the dressing map $d_\ell$.\\

\paragraph{Step 3~: Lower bound along the curve $\gamma$}~\\

Let us introduce a few more notations.
\begin{equation}\label{defrho1}
\rho(\tau) = \rho (\gamma(t,s,\tau))\;,\;  \rho =\rho(0), 
 \end{equation}
and 
\begin{equation}\label{defth}
\theta(\tau)=\theta_\tau =
( \frac{t(\tau)}{\rho(\tau)}, \frac{s(\tau)}{\rho^\ell(\tau)})
=(\widetilde t(\tau), \widetilde s(\tau))\;,
\end{equation}
with   $\theta =\theta(0)= (\widetilde t,\widetilde s)$.\\
We note that with the above notations
\begin{equation}
\varphi (t,s) = \rho^m \widetilde \varphi (\widetilde t,\widetilde
s)\;.
\end{equation}

We want to show (ii)$_{(c)}$ in Assumption \ref{Ass1} and first decompose
 the expression which has  to be  estimated from below by writing~:
\begin{equation}\label{decomp}
\rho(\tau)^m \widetilde  \varphi (\theta_\tau) - \rho^m
 \widetilde \varphi(\theta) = 
(\rho(\tau)^m - \rho^m) \widetilde \varphi (\theta)
 + \rho(\tau)^m (\widetilde \varphi(\theta_\tau) - \widetilde \varphi
 (\theta))\;,
\end{equation}
and will obtain a lower bound for each term of the r.h.s. of \eqref{decomp}.

\paragraph{Analysis of $\rho(\tau)^m - \rho^m$\,.}~\\
Using \eqref{inequt1} and 
$$
\rho(\tau)^m - \rho^m =
\left( \rho(\tau)^{2\ell} \right)^{\frac{m}{2\ell}}  - \left(
\rho^{2\ell} \right) ^{\frac{m}{2\ell}}
 \geq \left (\rho(\tau)^{2\ell} -\rho^{2\ell}\right)^{\frac{m}{2\ell}}
 \;,
$$
where we note that $m\geq 2 \ell$, 
we deduce
\begin{equation}\label{3.14a}
\rho(\tau)^m - \rho^m \geq 2^{ - m} \tau^m\;,\;\forall \tau \geq 0\;.
\end{equation}
Finally\footnote{ In the case, when $c=0$, then if $ds>0$, we have 
$$ \psi(\tau) = (s+d\tau)^2 -s^2 \geq \tau^2 \geq \tau^{2\ell}\;,\; \forall \tau \in [0,1].$$}, using  the lower bound for $\widetilde \varphi$, we get  for the first term of the r.h.s. of \eqref{decomp} the inequality~:
\begin{equation}\label{3.14b}
 (\rho(\tau)^m - \rho^m)
\widetilde \varphi (\theta) \geq 2^{-m} \mu \tau^m\;.
\end{equation}

\paragraph{Analysis of $\widetilde \varphi(\theta_\tau)-\widetilde
  \varphi (\theta)$\,.}~\\
Having in mind our assumptions on the variation of  $\widetilde \varphi$ on
$[\theta_1,\theta_2]$, we have simply 
 to prove here~:
\begin{itemize}
\item If $(\widetilde t,\widetilde s)\in [\theta_1,\widehat \theta_1]$,
  the function $\tau \mapsto \theta_\tau$ is non decreasing.
\item If $(\widetilde t,\widetilde s)\in [\widehat \theta_1,\theta_2]$, the
  function $\tau \mapsto \theta_\tau$  is non increasing.
\end{itemize}
But this is immediate after having sent the initial picture by the 
 dressing map \eqref{dressing}.

\section{Control in the case
 of the sectors $\mathcal V_j^+ \cup \mathcal R_{(a_j,b_j)}$}\label{s5a}
\paragraph{Provisory assumption}.\\
{\bf For the control, of the order of the zeroes, we assume for the
 moment that the $p$ appearing in \eqref{4.1} (Assumption
 \ref{H2}) satisfies \begin{equation}\label{p<m}
1\leq p\leq m\;.\end{equation}  This will be removed later 
(see Subsection \ref{sspm}).}\\

In comparison with the quasielliptic situation, the only point is that
 the condition of positivity could be not  satisfied at
 one or two ends. Note that we keep Assumption~\ref{Prop1}~iii).
The essential idea is then to improve the second part
 of the lower bound for  \eqref{decomp} using  the fifth
 item of Assumption \ref{H2} (i.e. the lower bound \eqref{4.1}).

Having in mind what we have done before, it remains to analyze
``rays'' coming from points $(t,s)$ which are close to $\mathcal R_{(a_j,b_j)}$.
The key estimate is that there exists $\mu_0 >0$  and a ``sector'' $\mathcal V$ 
neighboring  this ``end-ray'' such that
\begin{equation}\label{keyes}
\widetilde \varphi (\theta_\tau) - \widetilde \varphi (\theta)
 \geq \mu_0 \left( \frac{\tau}{\rho(\tau)}\right)^m\;,\; \forall
 (t,s)\in \mathcal V\;,\; \forall \tau \in [0,1]\;.
\end{equation}
But this proof is immediate from \eqref{4.1},  if we  control more
 precisely 
 the positivity of $\theta(\tau)-\theta(0)$, as we shall see in \eqref{319}.\\

\paragraph{Quantitative control of the positivity}~\\

We treat three typical cases.\\

\paragraph{Case 1 : $c>0, d>0$, $a_1>0, b_1 <0$}~\\
We assume that $\widetilde \varphi (a_1,b_1) =0$ and
 it is enough to control the trajectories starting
 from $(s,t)\in \sigma$ where the ``sector''  $\sigma$ is defined by the 
 condition that the  corresponding $(\widetilde t , \widetilde s )$ close to
 $(a_1,b_1)$ and between $(a_1,b_1)$ and $(c,d)$~:
\begin{equation}\label{defsigma}
(t,s)\in \sigma \mbox{\;iff\;} (\widetilde t,\widetilde s)\in \left(\theta_1, \theta_1
 +\epsilon_1\right) \,
,
\end{equation}
for some $\epsilon_1 >0$.\\
The point here is that one has a regular parametrization
 of the disto-circle by the $\tilde s$ variable and that we stay in
 the half-right plan.\\
By Formula \eqref{stildeprimea}, we observe that for $(t,s)$
 as before, we have
$$
\widetilde s' (\tau) \geq 0\;,\; \forall \tau \in [0,1]\;,
$$
Here we have used that  $\Delta_\ell(c,d,a_1,b_1)  < 0$. 
We also get
\begin{equation}
\widetilde \varphi(\theta_\tau)-\widetilde
  \varphi (\theta) \geq 0\;.
\end{equation}
We now use \eqref{stildeprimeb} and observe that 
 $\Delta_\ell(\widetilde t(\tau),\widetilde s(\tau), c,d)$ is increasing 
as $\tau$
 increases (together with $\widetilde s(\tau)$).
\\
Let us observe the trivial inequality 
\begin{equation}\label{minortriv}
1\geq \frac{t(\tau)}{\rho(\tau)}\geq c  \frac {\tau}{\rho(\tau)}\;.
\end{equation}
So there exists $\rho$ and for any $(s,t)\in \sigma$, $\tau(t,s)$
 such that $$ \Delta_\ell((\widetilde t(\tau),\widetilde
 s(\tau));(c,d))  
 \geq \rho\;,$$  for all $\tau \in [0,\tau(t,s)]$
 and $
 \Delta_\ell((\widetilde t(\tau),\widetilde s(\tau));(c,d)) 
 \leq \rho$ if $\tau \in [\tau(t,s),1]$.\\
In the first interval, we observe that $\widetilde t$ is semibounded
 on $\sigma$, so~:
\begin{equation}\label{minsa}
\widetilde s(\tau) - \widetilde s (0) \geq \mu_1 \frac{\tau}{\rho(\tau)}\;,\;
 \forall \tau \in [0,\tau(t,s)]\;,
\end{equation}
where $\mu_1= \ell (\inf_{(t,s)\in \sigma} \widetilde t^{2 \ell -1}) c^{1-\ell} \rho$.\\

This gives a uniform lower bound for $\tau \in [0,\tau(t,s)]$,
 but for $\tau \geq \tau(t,s)$ we have a  uniform lower bound
 (if we choose $\rho$ small enough) 
 of $\widetilde s (\tau) - \widetilde s(0)$ by a  constant $\alpha_\rho >0$.
 This together with \eqref{minortriv}, gives the existence of $\mu_2 >0 $ such that
\begin{equation}\label{minsb}
\widetilde s(\tau) - \widetilde s (0) \geq \mu_2 \frac{\tau}{\rho(\tau)}\;,\;
 \forall \tau \in [\tau(t,s), 1]\;.
\end{equation}
So we have finally shown that there exists $\mu_\sigma >0$
 such that 
\begin{equation}\label{minsc}
\widetilde s(\tau) - \widetilde s (0) \geq \mu_\sigma \frac{\tau}{\rho(\tau)}\;,\;
 \forall \tau \in [0, 1]\;.
\end{equation}

\paragraph{Case 2 :  $c>0, d>0$, $a_2<0, b_2 > 0$}~\\
We assume that $\widetilde \varphi (a_2,b_2) =0$
and
 it is enough to control the trajectories starting
 from $(s,t)$
 with corresponding $(\widetilde s (0), \widetilde t (0))$ close to
 $(a_2,b_2)$ and between $(c,d)$  and $(a_2,b_2)$ .\\
Here $\widetilde t (\tau)$ may change of sign along the trajectory
and we better parametrize the disto-circle by the variable $\widetilde
t$.\\
Here we use \eqref{ttildeprimea} for observing that $\widetilde
t(\tau)$
 is this time increasing (note that $\Delta_{\ell}( c,d,a_2,b_2) >0$
 and $s(\tau) \geq \varrho_0 >0$).

A similar argument to the one leading to \eqref{minsa} 
 gives the uniform control
 of $|\widetilde t (\tau) - \widetilde t(0)|$
 from below. Here we can no more use \eqref{minortriv}
 but will instead use \eqref{inequt1} which implies
 $$
\frac{\tau}{\rho(\tau)} \leq 2\;,
$$
without to assume the positivity of $t$.\\

There are limiting cases when we shall use both parametrizations
 but this does not create considerable troubles. Typically, let us
 consider the following case.
\paragraph{Case 3 : $c=0$, $d=-1$, $a_1<0$, $b_1<0$}~\\
We have just to use in this case vertical escaping ``rays''.\\
 
So in all the cases we obtain that,  if we start (inside 
 $\mathcal V_j^+$) of a ``subsector'' $\sigma $,  
 whose closure does not meet the ``ray'' 
 $\mathcal R_{(c,d)}$ then, there exists $\mu_\sigma >0$
 such that 
\begin{equation}\label{319}
|\theta (\tau) - \theta (0)| \geq \mu_\sigma \frac{\tau}{\rho(\tau)}\;,\; \forall \tau \in [0,1]\;.
\end{equation}

So we have finally  proved
\begin{lemma}\label{compl}~\\
Assuming that  $\widetilde
 \varphi$ satisfies in  $\mathcal V_j^+$ Property \ref{Prop1}
 and  the non degeneracy Assumption \eqref{4.1} for its  zero 
 possibly appearing at $(a_j,b_j)$, then   there exists  $\mu>0$ such that
\begin{equation}\label{minorcomp}
\rho(\tau)^m \varphi(\theta(\tau)) - \rho(0)^m \varphi(\theta(0))
 \geq \mu\; \tau^m\;,\; \forall \tau \in [0,1]\;,
\end{equation}
with $(\rho(0),\theta(0))$ corresponding to a point of $\mathcal V_j^+\cup \mathcal
R_{(a_j,b_j)}$.
\end{lemma}
\begin{remark}~\\
A similar result can be obtained for  $\mathcal V_j^+\cup \mathcal
R_{(a_{j+1},b_{j+1})}$ with a zero at $(a_{j+1}, b_{j+1})$.
\end{remark}
\section{Control for the ``sectors'' where $\widetilde \varphi$
 is non positive} \label{s5}
\subsection{Main case}
We consider first the case of small opening ``sectors'' $\mathcal V_{k}^-$.
More precisely, we assume that
\begin{equation}\label{5.1}
\langle (c,d)\,|\, (a_j,b_j)\rangle_{\ell} >0 \;,\; \mbox{ for } j=k,
k+1\;.
\end{equation}
\subsubsection{The subcase $c\neq 0$}
We start with the case when 
\begin{equation}\label{hypc}
c\neq 0\;,
\end{equation} 
We keep the same notations (change the labelling
by taking $k=1$)  but this time $(c,d)$ is a point of
$\mathcal S$ where $\widetilde \varphi$ is a minimum. The map $\gamma_-$
 is defined as follows.\\

If $(t,s)\in D(\frac 14)$ is such that $(\widetilde t,\widetilde s)\in
  \left((a_1,b_1), (c,d)\right)$, we let
\begin{enumerate}
\item
$\gamma_- =\gamma_1$ with 
\begin{equation}\label{gam1}\begin{array}{ll}
\gamma_1 (t,s,\tau) &:= (t(\tau),s(\tau))\\
&:= \left(t- c \tau, s +
\frac{d}{f_\ell(c)} (f_\ell(t(\tau) ) - f_\ell(t))\right)\;,\; \forall \tau \in
     [0,\tau_1(t,s)]\;,
\end{array}
\end{equation}
where $\tau_1(t,s)=\tau_1$ is the smallest $\tau$ such that $\gamma_1(t,s,\tau)
\in \mathcal R_{(a_1,b_1)}$. 
\item $\gamma_- =\gamma_2$ with 
\begin{equation}\label{gam2}
\gamma_2(t,s,\tau) = \gamma_+ (t(\tau_1), s(\tau_1), \tau
-\tau_1)\;,\;
\forall \tau \in [\tau_1,1]\;,
\end{equation}
where $ \gamma_+$ is the map constructed in \eqref{304a} in the
``sector'' $\mathcal V_{0}^+$ (the  ``sector'' preceding  $\mathcal V_1$ when turning anticlockwise).
\end{enumerate}

 \begin{figure}[h!]\label{picturesquare}
\begin{center}
  \includegraphics[width=8cm]{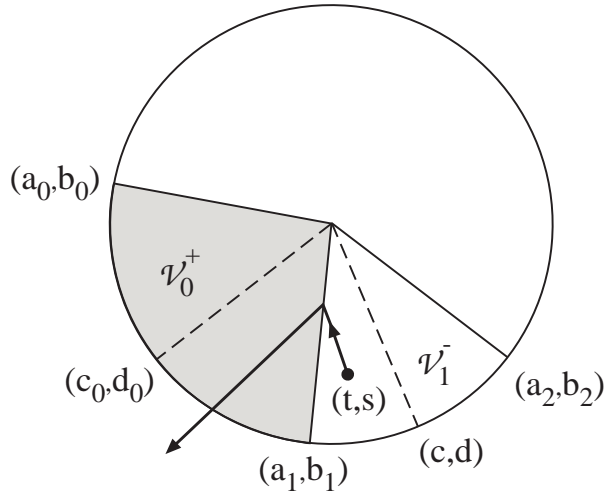}
\caption{The escaping ray (after application of the dressing map)}
  \end{center}
 \end{figure}

If $(\widetilde t,\widetilde s) \in \left((c,d),(a_2,b_2)\right)$, we
do the same construction
 with $\mathcal R_{(a_1,b_1)}$ replaced by $\mathcal R_{(a_2,b_2)}$.\\

Let us  compute $\tau_1(t,s)$. This
 $\tau_1$ is actually determined by writing that, at the corresponding
point
 $ 
(t(\tau_1),s(\tau_1))= (t_1(t,s),s_1(t,s))\;,
$
we should have
\begin{equation}\label{id1}
\Delta_\ell(a_1,b_1, t_1(t,s),s_1(t,s)) =0\;,
\end{equation}
and
\begin{equation}\label{id2}
\Delta_\ell (c,d, t_1(t,s),s_1(t,s)) = \Delta_\ell (c,d,t,s)\;.
\end{equation}
The first one \eqref{id1} expresses that we cross $\mathcal R_{a_1,b_1}$
 and the second one \eqref{id2} was observed in \eqref{inva}.
This leads to the determination of $\tau_1$
 by the formula
\begin{equation}
f_\ell(t-c \tau_1) = f_\ell( a_1)
\frac{\Delta_{\ell}(t,s,c,d)}{\Delta_{\ell}
(a_1,b_1,c,d)}\;.
\end{equation}
The uniqueness is obtained  by the monotonicity of $\tau \mapsto
 f_\ell(t-c \tau_1)$
 for $c\neq 0$
 and the existence is a consequence of the transversality of
 $\mathcal R_{a_1,b_1}$ and the disto-parallel to $\mathcal R_{c,d}$
 which is expressed by the condition
\begin{equation}
\Delta_\ell(a_1,b_1,c,d)\neq 0\;.
\end{equation}

{\bf It remains  to control the Jacobian of the map $\gamma$.}\\
It is immediate to see that the Jacobian of $(t,s)\mapsto \gamma (t,s,\tau)$
 is equal to one when $\tau < \tau_1 (t,s)$. Let us look at the more difficult case when 
$$\tau > \tau_1(t,s)\;. $$
Under this condition, we write
$$
\gamma_2(t,s,\tau) = (t_2(t,s,\tau), s_2(t,s,\tau)))\;,
$$
but omit the reference to $\tau$ in the computations below and prefer
to think of $t_1$, $t_2$ and $\tau_1$ as functions of $(t,s)$.\\
In order to compute the Jacobian, it is enough to compare the two
$2$-forms  
 $dt_2\wedge ds_2$ and $dt \wedge ds$.\\
We have in addition to \eqref{id1} and \eqref{id2}
\begin{equation}\label{identities}
\begin{array}{l}
t_1(t,s) = t - c \tau_1(t,s)\;,\\
t_2(t,s)  = t_1(t,s) + \hat c (\tau - \tau_1(t,s))\;,\\
s_2(t,s) = s_1(t,s) + \frac{\hat d}{f_\ell (\hat c)} (f_\ell(t_2(t,s)) - f_\ell(t_1(t,s)))\;.
\end{array}
\end{equation}
Here $(\hat c,\hat d)$ denotes the coordinates (and we assume for the
moment\footnote{See later for  what we have to do when $\hat c =0$.}  that $\hat c \neq 0$) of the maximum of
$\widetilde \varphi$ in the arc (corresponding to $\mathcal V_{k-1}^+$
 in the initial notations) preceding the arc $\left((a_1,b_1),(a_2,b_2)\right)$.\\

The two first identities imply first (by \eqref{id1}) that
\begin{equation}\label{id3}
dt_1\wedge ds_1 =0\;,
\end{equation}
and that (by \eqref{id1} and \eqref{id2}) one can express the
 $1$-forms $ds_1$ and $t_1(t,s)^{\ell -1} dt_1$
 as linear combinations of the $1$-forms  $ds$ and $f'_\ell(t) dt$. In particular we get~:
\begin{equation}\label{id4}
ds_1 = \frac{b_1}{f_\ell(a_1)} f'_\ell(t_1)\, dt_1\;,
\end{equation}
and
\begin{equation}\label{id5}
\frac{ 1}{ b_1} \Delta_\ell(c,d,a_1,b_1)\, ds_1
 = f_\ell(c) \,ds - d f'_\ell(t)\, dt\;.
\end{equation}
Let us start the computation of $dt_2\wedge ds_2$ using the two last identities.We obtain (using the rules of the exterior calculus)
$$
\begin{array}{ll}
dt_2\wedge ds_2 & =\left( dt_1 -\hat c d\tau_1\right)\wedge \left(ds_1
- \frac{\hat d}{f_\ell(\hat c)} d f_\ell (t_1(t,s)))\right)\\
& = - \hat c d \tau_1 \wedge 
\left(ds_1 - \frac{\hat d}{f_\ell (\hat c)}\, d (f_\ell(t_1(t,s)))\right)\\
& = \delta_1 \; d \tau_1 \wedge d s_1 \;, 
\end{array}
$$
with
$$
\delta_1= - |\hat c|^{1-\ell}\, b_1^{-1} \Delta_\ell (\hat c,\hat d,a_1,b_1)
\;.
$$
It remains to use the first relation of \eqref{identities}
 which gives~:
$$
d \tau_1= \frac 1c\, dt -\frac 1c\, dt_1\;,
$$
and
$$
dt_2\wedge ds_2= \frac {\delta_1}{c} \;  dt \wedge ds_1\;.
$$
Using \eqref{id5}, we finally obtain the following lemma.
\begin{lemma}\label{lbroken}~\\
If $\tau >\tau_1(t,s)$, then
\begin{equation}
dt_2\wedge ds_2 =\delta \; dt \wedge ds\;,
\end{equation}
with
\begin{equation}
\delta := - |\hat c|^{1-\ell} |c|^{\ell -1} \frac{ \Delta_\ell(\hat
  c,\hat d, a_1,b_1)}{\Delta_\ell( c,d,a_1,b_1)}\;.
\end{equation}
\end{lemma}
So the Jacobian is equal to $\delta$, hence constant,  and non
zero.
 The fact that $ \Delta_\ell(\hat
  c,\hat d, a_1,b_1)$ is not zero is the consequence  of the
  assumption on the zeros of $\widetilde \varphi$.

\begin{remark}~\\
The existence of this lower bound of the Jacobian  is probably the key
 point.  It is shown (see \cite{Tr}) in the analytic case
 that one can always find a $\gamma$ satisfying all the
 assumptions except this lower bound of the
 Jacobian by simply considering the flow
 associated with $\frac{1}{|\nabla \varphi|} \nabla \varphi$.
\end{remark}

Now if we observe that $(-\varphi)$ has in $\mathcal V_k^-$ the properties that
 $\varphi$ had in $\mathcal V_1^+$, we get rather easily   the existence of $\mu >0$ such that
\begin{equation}\label{5.5}
\rho(\tau)^m \widetilde \varphi (\theta_\tau) - \rho^m \widetilde
\varphi(\theta) 
 \geq \mu \tau^m\;, \forall \tau \in [0,\tau_1(t,s)]\;.
\end{equation}
Here we can indeed use Lemma \ref{compl} (after exchange of the roles  of $\theta_\tau$ and $\theta$).\\
We note  also that, for $\tau > \tau_1(t,s)$, we are in a region where
 $\varphi$ is positive so we can apply what we have done in this
 case. In particular, we obtain (see \eqref{3.14b}) with
 $\tau_1=\tau_1(t,s)$, 
\begin{equation}\label{5.6}
\rho(\tau)^m \widetilde \varphi(\theta_\tau) - \rho(\tau_1) \widetilde
 \varphi (\theta_{\tau_1}) \geq \mu (\tau - \tau_1)^m\;.
\end{equation}
One can also observe (see\footnote{ Here we use Assumption \eqref{5.1}.} \eqref{3.14a}) that
\begin{equation}\label{5.7}
\rho(\tau)^m - \rho^m \leq - (\frac \tau 2)^m\;,\; \forall \tau \in
    [0,\tau_1(t,s)]\;,
\end{equation}
which implies in particular the upperbound
\begin{equation}\label{5.7.a}
\tau_1(t,s) \leq 2\rho \leq \frac 12\;,
\end{equation}
and the inequality (see again \eqref{3.14a})
\begin{equation}\label{5.8}
\rho(\tau)^m - \rho(\tau_1)^m \geq (\frac{\tau - \tau_1}{2})^m\;,\; \forall \tau
\in [\tau_1,1]\;.
\end{equation}
~From \eqref{5.7.a} and \eqref{5.8}, one obtains 
that
\begin{equation}\label{5.9}
\gamma (t,s,1)\not\in D(\frac 14)\;.
\end{equation}
So the escaping condition (i) of $(H_+(\alpha))$ is satisfied with $\omega=D(\frac 14) $.\\

On the other hand, we get from \eqref{5.5} and \eqref{5.6}, the
estimate
\begin{equation}\label{5.10}
\rho(\tau)^m \widetilde \varphi (\theta_\tau) - \rho^m \widetilde
\varphi(\theta)
 \geq 2^{-m} \mu \tau^m\;,\; \forall \tau \in [0,1]\;.
\end{equation}

\subsubsection{The subcase $c=0$.}
For definiteness, we assume that $c=0$ and $d=-1$ and look
 at initial date in the ``sector'' attached to the arc
 $\left((a_1,b_1)\,,\,(0,-1)\right)$.
 The point is just to choose a point $(c',d')$ in $\mathcal S$,
 with $c'>0$ and sufficiently close to $(c,d)=(0,-1)$ in order to 
 keep the condition 
$$
\langle (a_1,b_1)\,|\, (c',d')\rangle >0\;.
$$
We can then keep the previous construction
 with $(c,d)$ replaced by $(c',d')$.

\subsection{Remaining case.}
In order to achieve the proof of Theorem \ref{theogen} (i.e. to prove
 that $(H_+(\alpha))$ is satisfied), we have to treat
 the case when \eqref{5.1} is no more satisfied. So  $(s,t)\mapsto \langle (c,d) \,|\, (t(\tau),
 s(\tau))\rangle $ may change of
 sign on $\mathcal V_k^-$. We have avoided this problem in the case of
 $\mathcal V_k^+$
 by dividing the ``sector'' in smaller ``sectors'',
 but this is no more possible when $\widetilde \varphi$ is negative. 
In this general case, one will be obliged
 to add  a broken line to the two previously defined arcs (ingoing and
 escaping) in order to leave $\omega$. We will see that it is always
 possible
 using a broken line made of at most five segments at the price to
 take $\omega$ smaller.\\
We start  from a point $(t,s)$ in the intersection of $\omega  :=
 D(R)$, ($R>0$ small  enough) with  a ``sector'' associated with the arc $\left((a_1,b_1)\,,\,
 (c,d)\right)$ and  
 we divide this ``sector'' into $N$ ``sectors'' $S^j$ ($j=1,\dots,N$) of disto-angle
 $<\frac \pi 2$.
 They are delimited by ``rays'' attached to the sequence $(\widehat
 \eta_j, \widehat \zeta_j)$ ($j=0,\dots, N)$ in $\mathcal S$ going {\bf
 clockwise}
 with   $(\widehat
 \eta_0, \widehat \zeta_0)=(c,d)$,  $(\widehat
 \eta_N, \widehat \zeta_N)=(a_1,b_1)$.\\
For the commodity of the notations below, we also introduce
 $(\widehat
 \eta_{-1}, \widehat \zeta_{-1})=(a_2,b_2)$ and  $(\widehat
 \eta_{N+1}, \widehat \zeta_{N+1})=(c_0,d_0)$.\\
 It is clear that we can always do that with $N\leq 4$ and that
 we have previously treated the case $N=1$.\\
We need also another sequence $(\eta_j, \zeta_j)$ in $\mathcal S$,
 which will determine the various directions of the
 broken line
 and  will satisfy
\begin{equation}\label{hypsuit}
\left\{
\begin{array}{l}
\eta_j \neq 0\;,\\
(\eta_j,\zeta_j) \in \left((\widehat \eta_{j-2}, \widehat\zeta_{j-2})\,,\,(\widehat \eta_{j-1}, \widehat\zeta_{j-1})\right)\;,\\
\Delta_\ell (\eta_j,\zeta_j,\widehat \eta_j,\widehat \zeta_j)   \neq
0\;,\\
\Delta_\ell (\eta_{j+1},\zeta_{j+1},\widehat \eta_j,\widehat
\zeta_j)\neq 0\;.
\end{array}
\right.
\end{equation}
This will be satisfied  by taking $(\eta_j,\zeta_j)$ very close (but
distinct
 except possibly for $N=1$ (if $c\neq 0$)) to 
$(\widehat \eta_{j-1}, \widehat \zeta_{j-1})$ for $j=1,\dots,N$.\\
So we can always assure the property that starting form a point in the
disk
 on $\mathcal R_{\widehat \eta_j,\widehat \zeta_j}$ the ``straight''
 line parallel to $(\eta_{j+1},\zeta_{j+1})$
 will meet the ``ray'' $\mathcal R_{\widehat \eta_{j+1},\widehat
   \zeta_{j+1}}$ inside the ``disk'', for $j=1,\dots, N-1$.\\

We now explain how we construct the broken line starting
 from
 a point $(t,s)$ belonging to the first ``sector'' $S_1:=S(\widehat
 \eta_0,\widehat \zeta_0,\widehat \eta_1, \widehat \zeta_1)$. The
 other cases are simpler.\\

Starting from $(t,s)$ we follow for $\tau\geq 0$
 the parametrized ``line'' parallel to $(\eta_1,\zeta_1)$ till we meet
 at the time $\tau_1(t,s)$ the ``ray'' $\mathcal R_{\widehat
 \eta_1,\widehat \zeta_1}$ at
 the point $(t_1(t,s),s_1(t,s))$.\\
\\
If $N>1$ and  starting now  from $(t_1(t,s),t_2(t,s))$ we follow for $\tau\geq \tau_1(t,s)$
 the parametrized ``line'' parallel to $(\eta_2,\zeta_2)$ till we meet
 at the time $\tau_2(t,s)$ the ``ray'' $\mathcal R_{\widehat \eta_2,
 \widehat \zeta_2}$ at
 the point $(t_2(t,s),s_2(t,s)$.\\

By recursion, we arrive at the point $(t_N(t,s), s_N(t,s))$ on the ``ray''
$R_{\widehat \eta_N,\widehat \zeta_N }$ at time $\tau_N(t,s)$. We can
then escape along $(\eta_{N+1},\zeta_{N+1})$.\\

For generalizing, what we have done for $N=1$, we have now mainly to
 verify
 the following points~:
\begin{enumerate}
\item The condition of the Jacobians
 for $\tau \in ]\tau_{j}(t,s),\tau_{j+1}(t,s)[$~:\\  This will be
 analyzed in the appendix,
\item The control of the escape time~:\\ We have   to determine
a condition on $\omega$ 
 such that 
\begin{equation}
\tau_{j+1}(t,s) -\tau_{j} (t,s) \leq \frac{1}{4N}\;,
\end{equation}
\item The control of the dynamics~:\\
We have  to control
$\rho(\tau_{j}(t,s))$ for $j=1,\dots, N$, under the suitable
assumption that $(s,t)\in D(0,R_N)$.
\item The positivity of $\widetilde \varphi (\gamma(t,s,\tau)) -
  \widetilde \varphi (t,s)$ along the trajectory~:\\ We do not meet here new
  problems.
\end{enumerate}

For the proof of (ii) and (iii)  we observe that
\begin{equation}
\rho(\tau)^m + \left( \frac{\tau - \tau_i}{2}\right)^m \leq
\rho(\tau_i)^m\;,\; \forall \tau \in ]\tau_i,\tau_{i+1}]\;.
\end{equation}
So in particular
\begin{equation}
\rho(\tau_{i+1})^m + \left( \frac{\tau_{i+1} - \tau_i}{2}\right)^m \leq
\rho(\tau_i)^m\;.
\end{equation}
This gives 
\begin{equation}
\rho(\tau_{i+1}) \leq \rho(\tau_i)\;,
\end{equation}
and 
\begin{equation}
\tau_{i+1} - \tau_i \leq 2 \rho(\tau_i)\leq 2 \rho(0)\;.
\end{equation}
In addition, we have (see \eqref{5.8} for the case $N=1$)
\begin{equation}
\rho(\tau)^m - \rho(\tau_{N}(t,s))^m \geq2^{-m}  (\tau- \tau_{N}(t,s))^m\;.
\end{equation}
~From this we deduce (without to look for optimality) that $$R_N = \frac
18 \frac{1}{N}$$
 is enough for showing that the parametrized broken line has left
 $\omega = D(0,R_N)$ at time $\tau =1$.

\subsection{The case when $p> m$}\label{sspm}
Till now we have proved the main theorem under the additional
condition \eqref{p<m}. Without this additional condition, the
subellipticity which is obtained 
 is now $\inf(\frac 1m,\frac 1p)$ instead of $\frac 1m$.
The only change is that we get only the existence of $\mu >0$
 such that, for $\tau \in [0,1]$,  
\begin{equation}
\rho(\tau)^m \varphi(\theta_\tau) - \rho^m \varphi(\theta)
 \geq  \mu\, \inf\left( \tau^m, \rho(\tau)^m
 \left(\frac{\tau}{\rho(\tau)}
\right)^p\right)\;.
\end{equation}
But we have shown the existence of a constant $C_{\varphi}>0$ and the
existence of some open neighborhood of $0$ $\omega_{\varphi}$ such
that $$0\leq \tau \leq C_{\varphi} \, \rho(\tau)\;.$$
So we get
\begin{equation}
\rho(\tau)^m \varphi(\theta_\tau) - \rho^m \varphi(\theta)
 \geq  \mu\, \tau^{\sup(m,p)}\;.
\end{equation}

\section{The analytic case and $\ell\in \mathbb Q $ }\label{s6}
\subsection{The main theorem in the analytic case}
We keep the previous assumptions but now assume that
\begin{equation}
\ell =\ell_2/\ell_1\;,
\end{equation}
with $\ell_1$ and $\ell_2$ mutually prime integers, 
 and that $\varphi$ is analytic. 
 In this
case, the quasihomogeneity Assumption \eqref{quash} on $\varphi$ implies that $\varphi$ 
 is actually a polynomial and we can write $\varphi$ in the form
\begin{equation}\label{qhan}
\varphi(t,s) =\sum_{\ell_1 j+\ell_2 k =\ell_1 m
} a_{j,k} t^j s^k\;,
\end{equation}
where $(j,k)$ are integers and the $a_{j,k}$ are real.\\
We  can of course apply the  main theorem but 
 it is  nicer to have a criterion involving
 more directly the assumptions on $\varphi$ or on 
its restriction $\widehat \varphi$  of $\varphi$
 to the quasi-circle $$ \mathcal S_{\ell_1,\ell_2} :=\{
t^{2\ell_2} + s^{2\ell_1} =1\}\;.
$$ 
instead of the disto-circle $\mathcal S$.
There are absolutely no problems
 if the critical points or zeroes of $\widetilde \varphi$
 avoid $\{t=0\} \cup \{s=0\}$ but one should be
 more careful in order to analyze Condition \eqref{4.1} at the
 remaining points.\\

Let us show how this works in this case. We parametrize (this is a
$C^2$ parametrization) on $\mathcal S$ by $t$ and 
assume that we are close to $(0,1)$ for definiteness and
$\widetilde \varphi$ becomes locally with this parametrization the function
$$
t \mapsto \kappa(t)=\varphi (t, \sqrt{1 -t^{2\ell}})\;,$$
and we assume that
$\kappa(0)=0$ and that $\kappa$ is not identically $0$.
Suppose that we are on the side $\{t>0\}$. Then 
$$
\kappa(t)=\chi (t^{\frac {1}{\ell_1}})\;,
$$
where $\chi$ is the non identically zero  analytic function
$$
v\mapsto \chi(v)= \varphi ( v^{\ell_1}, \sqrt{1 - v^{2 \ell_2}})\;,
$$
with $\chi(0)=0$.\\
Now there exists $p>0$ such that $\chi^{(p)} (0)\neq 0$ and we get
 the existence of $C>0$ such that, in the neighborhood of $0$
$$
\forall v \geq 0, \forall v' \geq 0, |\chi(v)-\chi(v')| \geq \frac 1C
|v-v'|^p\;.
$$
Coming back to $\kappa $, we get for a constant $\widehat C>0$
$$
\forall t \geq 0, \forall t' \geq 0, |\kappa (t)-\kappa (t')| \geq \frac 1C
|t^{\frac{1}{\ell_1}}- (t')^{\frac{1}{\ell_1}}|^p \geq \frac{1}{\widehat C} 
 |t-t'|^p\;.
$$
So we have obtained the proof of  \eqref{4.1} for some\footnote{Using
  that $\varphi$ is a polynomial, one could get more information on
  $p$ if needed.} $p$.

\begin{theorem}\label{th6.1gen}~\\
Let $\varphi$ be a real analytic non identically $0$ quasihomogeneous function 
  satisfying  \eqref{condl}, \eqref{condlm},  and \eqref{quash},
 with $\ell=\ell_2/\ell_1$. If $ \varphi$ is
 strictly positive
 on $\mathbb R^2\setminus \{0\}$, then $\varphi$ satisfies $(H_+(\alpha))$ with $\alpha=\frac {\ell_1}{m}$
 and the system \eqref{sys} is $\frac {\ell_1}{m}$-microlocally subelliptic in
 $\{ \xi >0\}$.
\end{theorem}

\begin{theorem}\label{th6.2gen}~\\
Let $\varphi$ be a real analytic non identically $0$ quasihomogeneous function 
  satisfying  \eqref{condl}, \eqref{condlm},  and \eqref{quash}, with $\ell=\ell_2/\ell_1$.  
Suppose that $\varphi$ is not a negative function. Suppose in addition that~:\\
If $\mathcal S^-_k = [\theta_k,\theta_{k+1}]$ is a maximal arc  where
$\widehat \varphi$ is negative, then $\widehat \varphi'$ has a
  unique
 zero on $]\theta_k,\theta_{k+1}[$.\\
Then $\varphi$ satisfies $(H_+(\alpha))$ with $\alpha >0$. Hence the
system \eqref{sys} is microlocally subelliptic in $\{ \xi >0\}$.
\end{theorem}

\begin{example}~\\
We recover some examples treated by H.~Maire \cite{Mai4}
$$
\varphi (t,s)= t (s^2- t^{2\ell})\;,\; \ell \geq 1\;.
$$
Here $m=2\ell +1$, $p=1$ and we get the subellipticity with 
$\alpha =\frac{1}{2\ell +1}$. As observed in \cite{HeNi}, this result
is optimal and  the associated system
 is not maximally hypoelliptic when $\ell >1$. The maximal
 hypoellipticity
 would indeed imply $\alpha =\frac 13$.
\end{example}

\subsection{Around Journ\'e-Tr\'epreau's examples}
For
$$
\varphi(t,s) = - t^{2m} - t^2 s^{2p} + s^q\;,
$$
with
$$
m\geq 1\;,\; p\geq 2\;,\; q\geq \frac  {2mp}{m-1}\;,
$$
J.L. Journ\'e and J.M. Tr\'epreau show that, although the
Maire-Tr\`eves condition is  satisfied,  one cannot
 obtain a better $\rho$-subellipticity than
$$
\rho \leq - ( 1-\frac{2p}{q} - \frac 1m ) \frac{n-1}{4} + \frac{1}{2q}
+ \frac{m-1}{4mp}\;.
$$
The right hand side can become strictly negative. But if we impose 
the  
quasihomogeneity condition \eqref{qhan}, we get as a necessary condition~:
$$
( 1-\frac{2p}{q} - \frac 1m ) =0\;,
$$
which cancels the only possible negative term.\\

Indeed inside  this class ($m=2$, $p=2$), the authors can obtain the optimal
subellipticity for the example
$$
\varphi (t,s)= - t^4 - t^2 s^4 + s^q\;,
$$
with $q \geq 8$.\\
The optimal subellipticity is $\rho_q =\frac{3}{2q} - \frac{1}{16}$.
Here let us observe that the only quasihomogeneous case corresponds to
$q=8$
 and that in this case their result is coherent with our result. This
 example shows also
 that we loose the ``positive'' subellipticity for 
$q\geq 24$.
\subsection{Final remarks}

In the analytic case, the criterion of microhypoellipticity (proven by
Maire) at say $(0,0)$
 in the direction $\{\xi >0\}$ is that $\varphi$ has no local maximum
 in a neighborhood of $0$. When $\varphi$ is quasihomogeneous, we
 immediately
 see that, at a local maximum,  $\widetilde\varphi$ should be zero. So we should
 avoid
 the following situations~:
\begin{itemize}
\item $\widetilde \varphi \leq 0$, with
  $\widetilde\varphi^{-1}(]-\infty,0[)\neq \emptyset$.
\item $\widetilde \varphi$ has a local maximum equal to zero 
 on $\mathcal S$.
\end{itemize}
One can indeed verify that our assumptions exclude these two cases but
are unfortunately more restrictive.\\

Note also that it would be interesting the case when $\widetilde
\varphi$
 or $\widehat \varphi$
 have a strictly negative local maximum. We have not been able by our
 method to construct escaping curves satisfying all the conditions
 of the criterion in this case.

Finally, let us recall that the maximal hypoellipticity of these
systems was analyzed in \cite{HeNo, No1,No4,No5} and more recently in \cite{HeNi}.

\newpage 
\appendix

\section{Control of the Jacobian for a $(N+1)$-broken line}
We just adapt the proof done for a broken line with two pieces.  
We treat for definiteness the case $N=2$ and keep the same notations
as in the main text.

Starting from $(t,s)$ we consider
 the parametrized (disto)-line
$$
t_1(\tau,t,s) = t - \eta_1 \tau\;,\; s_1(\tau,t,s) =
 s + \frac{\zeta_1}{f_\ell(\eta_1)}
 \left(f_\ell(t_1(\tau,t,s)) - f_\ell(t)\right)\;.
$$
Here $(\eta_1,\zeta_1)$ is a point of $\mathcal S$ with $\eta_1 \neq
0$
 and, in our example (if in addition $c\neq 0$),  $(\eta_1,\zeta_1) = (c,d)$.\\
We refer to \eqref{hypsuit} for  the conditions on the  sequence
$(\eta_j,\zeta_j)$
 and $(\widehat \eta_j,\widehat \zeta_j)$. We simply add the equations
 which are satisfied at the singular points of the broken line. The
 intersection with the ``ray'' $\mathcal R_{\widehat \eta_1,\widehat \zeta_1}$ occurs at time $\tau_1(t,s)$
 and we get the first family of relations for the intersection point
 denoted by $(t_1(t,s), s_1(t,s))$.\\
\begin{equation}\label{E3}
\begin{array}{l}
\Delta_\ell (\widehat \eta_1,\widehat \zeta_1,
t_1(t,s),s_1(t,s))=0\;;\\
\Delta_\ell ( \eta_1,\zeta_1, t_1(t,s),s_1(t,s))= \Delta_\ell
(\eta_1,\zeta_1, t,s)\;;\\
t_1(t,s) = t - \eta_1 \tau_1(t,s)\;.
\end{array}
\end{equation}

Starting from $(t_1(t,s), s_1(t,s))$ we construct a new ingoing
``straight''\footnote{more precisely, straight after application of
  the dressing map,} half-line
 $\mathcal L_2(t,s)$
 disto-parallel to the disto-unit-vector $(\eta_2,\zeta_2)$.

 The straight line meets   $\mathcal R_{\widehat \eta_2,\widehat
   \zeta_2}$ at time $\tau_2(t,s)$  
and at the point $(t_2(t,s), s_2(t,s))$. So we get the second family of relations~:

\begin{equation}\label{E4}
\begin{array}{l}
\Delta_\ell (\widehat \eta_2,\widehat \zeta_2,
t_2(t,s),s_2(t,s))=0\;;\\
\Delta_\ell ( \eta_2,\zeta_2, t_2(t,s),s_2(t,s))= \Delta_\ell
(\eta_2,\zeta_2, t_1(t,s),s_1(t,s))\;;\\
t_2(t,s) = t_1(t,s) - \eta_2 (\tau_2(t,s)-\tau_1(t,s))\;.
\end{array}
\end{equation}

For $\tau >\tau_2(t,s)$, we now consider the last half-line, 
which is now assumed to be escaping  and starting from $(t_2(t,s), s_2(t,s))$
 and parallel to a new vector $(\eta_3,\zeta_3)$. At time $\tau$, we will be at the point
\begin{equation}\label{E5}
\begin{array}{ll}
t_3(\tau,t,s) &= t_2(t,s) + \eta_3 (\tau - \tau_2(t,s))\;,\\
 s_3(\tau,t,s)& = s_2(t,s) + \frac{\zeta_3}{f_\ell(\eta_3)} (f_\ell(t_3(\tau,t,s))
 - f_\ell (t_2))\;.
\end{array}
\end{equation}
It remains to control the Jacobian for the different values of
 $\tau\in [0,1]$, which can be done
 by the computation of the coefficient of the 2-form $d t_3 \wedge ds_3$
 on the $2$-form $dt\wedge ds$. We will see that this coefficient
 is locally constant.\\
The proof is similar to the proof of Lemma \ref{lbroken}.\\
We are actually able to give explicit formulas
 of these Jacobians, once the two sequences  $(\eta_j,\zeta_j)$ and $(\widehat \eta_j,\widehat
 \zeta_j)$ are fixed. \\
Let us treat for definiteness the case $N=2$ and look first
 at what is going on for $\tau \geq \tau_2(t,s)$.\\
We will show~:
\begin{equation}\label{E6}
dt_3\wedge ds_3 = \underline \delta_{32} \,d \tau_2\wedge ds_2
 =\underline \delta_{32}  \widehat \delta_{21}\, d \tau_1\wedge ds_1
 = \delta_{3,0}\, dt \wedge ds\;.
\end{equation}
Let us show the existence of $\underline \delta_3$. Using \eqref{E5},
 we first get
\begin{equation}\label{E7}
dt_3\wedge ds_3 = dt_3 \wedge \left( ds_2 -
\frac{\zeta_3}{f_\ell(\eta_3)} f'_\ell(t_2) dt_2\right)\;.
\end{equation}
Now, the first line of \eqref{E4}, gives that
\begin{equation}\label{E8}
 f_\ell(\widehat \eta_2) ds_2 - \widehat \zeta_2  f'_\ell(t_2) dt_2
= 0 \;.
\end{equation}
We then obtain easily
\begin{equation}\label{E9}
\underline \delta_{32} = - \eta_3 \left(1 - \frac{ f_\ell(\widehat
  \eta_2)\zeta_3}{f_\ell(\eta_3 )\widehat \zeta_2}\right)
= - |\eta_3|^{1-\ell} {\widehat \zeta_2}^{\,-1} \Delta_\ell
  (\eta_3,\zeta_3,\widehat \eta_2, \widehat \zeta_2)\;, 
\end{equation}
which should be non zero in our construction.\\
The second step is to show that
\begin{equation}
d \tau_2 \wedge ds_2 = \widehat \delta_{21}\, d\tau_1 \wedge ds_1\;.
\end{equation}
The differentiation of the third line of \eqref{E4} gives~:
\begin{equation}
d \tau_2 \wedge ds_2 = (\frac{1}{\eta_2} dt_1 + d \tau_1) \wedge
ds_2\;.
\end{equation}
The differentiation of the second  line of \eqref{E4} (together with
 \eqref{E8} ) permits us  to express
 $ds_2$ as a combination of $ds_1$ and $f'_\ell(t_1) dt_1$.\\
We get indeed
$$
ds_2  = \frac{\widehat \zeta_2}{\Delta_\ell (\eta_2,\zeta_2,\widehat
  \eta_2,\widehat \zeta_2)}  \left(f_\ell(\eta_2) ds_1
-\zeta_2 f'_\ell( t_1) dt_1\right)\;,
$$
and  by differentiating the first line of \eqref{E3}, we get
\begin{equation}
\widehat \delta_{21} =  \frac{\widehat \zeta_2}{\Delta_\ell (\eta_2,\zeta_2,\widehat
  \eta_2,\widehat \zeta_2)} \frac{\Delta_{\ell}(\eta_2,\zeta_2, \widehat 
  \eta_1,\widehat \zeta_1)}{\widehat \zeta_1}\;.
\end{equation}
The last step is to show
$$
d\tau_1\wedge ds_1 = \overline  \delta_{10}\, dt\wedge ds\;.
$$
The differentiation of the first line of \eqref{E3} leads to
$$
d\tau_1\wedge ds_1 = \frac{1}{\eta_1} \, dt\wedge ds_1\;.
$$
It remains then to use the two first lines of \eqref{E3} and
we obtain
\begin{equation}
\overline \delta_{1,0} = |\eta_1|^{\ell-1}  \frac
{\widehat \zeta_1}{\Delta_\ell(\eta_1,\zeta_1,\widehat \eta_1, \widehat \zeta_1)}\;.
\end{equation}
So finally, we have obtained, observing that 
$$
\delta_{3,0} =\underline \delta_{3,2} \,\widehat \delta_{2,1}\, \overline \delta_{1,0}\,, $$
and 
consequently
\begin{equation}
\delta_{3,0}= -
|\eta_3|^{1-\ell}|\eta_1|^{\ell -1}
 \frac{ \Delta_\ell
  (\eta_3,\zeta_3,\widehat \eta_2, \widehat \zeta_2) }{\Delta_\ell (\eta_2,\zeta_2,\widehat
  \eta_2,\widehat \zeta_2)}  \frac{\Delta_{\ell}(\eta_2,\zeta_2, \widehat 
  \eta_1,\widehat \zeta_1)}{\Delta_\ell (\eta_1,\zeta_1,\widehat
  \eta_1,\widehat \zeta_1)}\;.
\end{equation}
Let us recapitulate what we have obtained  in the case of the three-broken line.
\begin{enumerate}
\item
For $\tau \in [0,\tau_1(t,s)[$, the Jacobian is one~:
\begin{equation}
\delta^{(1)} = \delta_{1,0} =1\;.
\end{equation}
\item 
 For $\tau \in
]\tau_1(t,s),\tau_2(t,s)[$, the Jacobian (see the main text) is 
\begin{equation}
\delta^{(2)} := \delta_{2,0}= - |\eta_2|^{1-\ell}\, |\eta_1|^{\ell -1} 
 \frac{\Delta_{\ell}(\eta_2,\zeta_2, \widehat 
  \eta_1,\widehat \zeta_1)}{\Delta_\ell (\eta_1,\zeta_1,\widehat
  \eta_1,\widehat \zeta_1)}\;.
\end{equation}
\item 
For $\tau \in ]\tau_2(t,s), +\infty[$, the Jacobian is given by
$\delta^{(3)}= \delta_{3,0}$.
\end{enumerate}
\paragraph{Generalization}~\\
It is now not too difficult to extend the formula in the case of  $N$
reflexions.
\begin{equation}
\delta^{(N+1)} = - |\eta_{N+1}|^{1-\ell} |\eta_1|^{\ell
  -1}\,\Pi_{j=1}^{N}\left( \frac{\Delta_{\ell} (\eta_{j+1}, \zeta_{j+1},
    \widehat \eta_j, \widehat \zeta_j)}{\Delta_{\ell} (\eta_{j}, \zeta_{j},
    \widehat \eta_j, \widehat \zeta_j)}\right)\;.
\end{equation}

\bibliographystyle{plain}

\end{document}